\documentclass[11pt]{article}

\usepackage[latin1]{inputenc}
\usepackage{amsmath,amssymb}
\usepackage{latexsym}

\newtheorem{theorem}{Theorem}[section]
\newtheorem{corollary}[theorem]{Corollary}
\newtheorem{lemma}[theorem]{Lemma}
\newtheorem{proposition}[theorem]{Proposition}

\newtheorem{remark}[theorem]{Remark}
\numberwithin{equation}{section}

\def\proof{{\medskip\noindent {\bf Proof. }}}

\def\qed{{\hfill $\square$ \bigskip}}

\def\square{{\vcenter{\vbox{\hrule height.3pt
        \hbox{\vrule width.3pt height5pt \kern5pt
           \vrule width.3pt}
        \hrule height.3pt}}}}

 \def\sB {{\cal B}} 
  \def\sF {{\cal F}}
  
  \def\sL {{\cal L}}

\def\wt{\widetilde}
\def\wh{\widehat}
\def\ol{\overline}
\def\E{{\mathbb E}}
\def\P{{\mathbb P}}
\def\norm#1{{\Vert #1 \Vert}}

\def\lam{{\lambda}}

\def\bea{\begin{align*}}
\def\eea{\end{align*}}
\def\bee{\begin{equation}}
\def\eee{\end{equation}}

\def\R{{\mathbb R}}
\def\E{{{\mathbb E}\,}}
\def\P{{\mathbb P}}

\def\lam{{\lambda}}

\def\al{{\alpha}}
\def\grad{{\nabla}}
\def\proof{{\medskip\noindent {\bf Proof. }}}

\def\qed{{\hfill $\square$ \bigskip}}

\def\eps{\varepsilon}
\def\vp{\varphi}

\def\norm#1{\Vert #1 \Vert}

 \def\qq {\qquad}

\def\wt{\widetilde}
\def\ol{\overline}

\def\wh{\widehat}

\def\bs{\bigskip}

\def\square{{\vcenter{\vbox{\hrule height.3pt
        \hbox{\vrule width.3pt height5pt \kern5pt
           \vrule width.3pt}
        \hrule height.3pt}}}}

\def\tfrac#1#2{{\textstyle {\frac{#1}{#2}}}}

\def\tlint{{- \kern-0.85em \int \kern-0.2em}}  
\def\dlint{{- \kern-1.05em \int \kern-0.4em}}  

 \def\sB {{\cal B}} 
  \def\sF {{\cal F}}
  
  \def\sL {{\cal L}}

\def\nn{{\nonumber}}

\begin{document}

%
%
%
%
%
%


\thispagestyle{empty}
\vspace*{0.7cm}
\centerline{\Large\bf Regularity of Harmonic Functions}
\medskip
\centerline{\Large\bf for a Class of Singular Stable-like Processes}

\vskip 0.3truein
\centerline{\large {\bf Richard F. Bass}\footnote{Research partially supported by NSF grant DMS-0601783.}
 \qquad and \qquad {\bf Zhen-Qing Chen}\footnote{Research partially supported by NSF
grant DMS-0600206.} }

\vskip 0.2truein

 { \halign{{}\hskip 1truein  # \hfill\hskip 0.6truein & # \hfill
  \cr
 Department of Mathematics   & Department of Mathematics   \cr
 University of Connecticut   &  University of Washington \cr
 Storrs, CT 06269-3009, USA  &  Seattle, WA 98195, USA \cr
 {\texttt bass@math.uconn.edu} &  {\texttt
 zchen@math.washington.edu} \cr}}

 \vskip 0.2truein

 \centerline{(April 21, 2009)}
\vskip 0.2truein

\begin{abstract}
\noindent
We consider the system  of stochastic differential equations $$dX_t=A(X_{t-})\, dZ_t,$$
where $Z_t^1, \ldots, Z^d_t$ are independent one-dimension\-al symmetric stable processes of
order $\al$, and the matrix-valued function $A$ is
 bounded,
 continuous and everywhere non-degenerate. We show that
  bounded
harmonic functions associated with $X$ are
H\"older continuous, but a Harnack inequality need not hold. The
L\'evy measure associated with the vector-valued process $Z$ is
highly singular.

\bigskip
\medskip


 \bigskip

\end{abstract}


\section{Introduction}

A one-dimensional symmetric stable process of index $\al\in (0,2)$
is the L\'evy process taking values in $\R$ with no drift, no
Gaussian part, and L\'evy measure
$$n(dh)=c_1/|h|^{1+\al}\, dh.$$
Let $Z_t=(Z^1_t, \ldots, Z^d_t)$ be a vector of $d$ independent
one-dimension\-al symmetric stable processes of  index $\al$.
Consider the system of stochastic differential equations
\bee\label{RE1}
 dX^i_t=\sum_{j=1}^d A_{ij}(X_{t-})\, dZ^j_t, \qq
X^i_0=x^i_0, \qq i=1, \ldots, d, \eee where $x_0=(x^1_0, \ldots,
x^d_0)\in\R^d$ and $A(x)$ is
 a bounded
 $d\times d$ matrix-valued function
that is continuous in $x$ and everywhere non-degenerate, that is,
  the determinant ${\rm det} (A(x))\not= 0$
  for all $x$. The main
result of \cite{system} is that under these conditions there is a
unique weak solution to the system \eqref{RE1}
 and the family $\{X, \P^{x_0}, x_0\in\R^d\}$ forms a
 strong Markov process on $\R^d$.
 The process $X$ may
be referred to as stable-like because it possesses an approximate
scaling property similar to the stable processes; see \cite{tang}
and \cite{CK03} for other examples where the term stable-like has
been used. The system \eqref{RE1} has been suggested as a possible
model for a financial market with jumps in the security prices
(\cite{CS}).
Note that by Proposition 4.1 of \cite{system}, the infinitesimal
generator of the Markov process $X$ determined by \eqref{RE1} is
 \begin{equation}\label{e:generator}
 \sL f(x)=\sum_{j=1}^d
\int_{\R\setminus\{0\}} \left( f(x+a_j(x)w)-f(x)-w1_{\{|w|\leq 1 \}}
\grad f(x)\cdot a_j(x)\right) \frac{c_1}{|w|^{1+\al}} dw,
 \end{equation}
where $a_j(x)$ is the $j^{{th}}$ column of the matrix $A (x)$.
Associated with the operator $\sL$ is the symbol
$$\ell(x,u):= c_2\sum_{j=1}^d |u\cdot a_j(x)|^\al,
\qquad x, u \in \R^d.
$$
This means
$$\sL f(x)= \int_{\R^d}  \ell(x,u) e^{-iu\cdot x} \wh f(u)\, du,$$
where $\wh f$ denotes the Fourier transform of $f$. This is an
example of a pseudodifferential operator with
singular
 state-dependent symbol.

We say that a function $h$ that is bounded in $\R^d$ is harmonic (with respect
to $X$) in a
domain $D$ if $h(X_{t\land \tau_D})$ is a martingale with respect to
$\P^x$ for every $x\in D$, where $\tau_D$ is the time of first
exit from $D$.
The process $X$ is shown to have no explosions in finite time  in \cite{system} and
when $D$ is   bounded, it is easy to see from \eqref{RE1} that
$\P^x(\tau_D<\infty)=1$ for every $x\in D$. So by the bounded
convergence theorem and the strong Markov property of $X$, a bounded
function $h$ on $\R^d$  is harmonic in a bounded domain $D$ if
and only if
$$ h(x)=\E^x[ h(X_{\tau_D })] \qquad \hbox{for every } x\in D.
$$
Consequently, every bounded harmonic function in a bounded domain
$D$ is the difference of two non-negative bounded harmonic functions
in $D$.
It follows from Proposition 4.1 of \cite{system} that a bounded
$C^2$ function $u$ is harmonic in $D$ if and only if $\sL u=0$ in
$D$.

 The main goal of this paper is to prove
 the H\"older continuity of functions which are bounded and harmonic with
respect to $X$ in a domain.

There are two reasons why the H\"older continuity is perhaps a bit unexpected.
Consider the case where $A$ is identically equal to the identity matrix, and
so $X\equiv Z$. Even in this case a Harnack inequality may fail; see
Section \ref{S4}. Nevertheless the H\"older continuity of the harmonic
functions holds. The other reason is that the process $Z$ is quite
singular. It is a L\'evy process, but the support of its L\'evy measure
is the union of the coordinate axes. By contrast, the support of the L\'evy
measure for a $d$-dimensional
(rotationally)
 symmetric stable process is all of $\R^d$, a much
more tractable situation.

The key to our method is the technique of Krylov-Safonov as given, for
example, in the exposition in \cite{Ba97}. The most difficult
step in our proof is the proof of a support theorem for $X$; this is given
in Section \ref{S3}. We remark that the current paper is the first one
where the full strength of the Krylov-Safonov technique has been used in
the context of pure jump processes.

 For a Borel subset $C\subset \R^d$,
let $T_C:=\inf\{t\geq 0: X_t\in C\}$ and $\tau_C:=\inf\{t\geq 0:
X_t\notin C\}$ be the first entrance and departure time of $C$ by
$X$. Let $|C|$ denote the Lebesgue measure of a Borel set $C$.
 The open ball of radius $r$ centered at $x$ will be denoted as $B(x,r)$.
 The paths of $Z_t$ are right continuous with left limits. We write
$$Z_{t-}:=\lim_{s\uparrow t, s<t} Z_s, \qq \Delta Z_t:=Z_t-Z_{t-},$$
and similarly $X_{t-}$ and $\Delta X_t$.
The letter $c$ with a subscript denotes  a positive finite constant whose exact value
is unimportant and may vary from one usage to the next.
Constant $c$ typically  depends on $\alpha$ and $d$, but for
convenience this dependence will not be explicitly mentioned
throughout the paper.

\section{Regularity}\label{S3}

For $1\leq i\leq d$, denote by $e_i$ the unit vector in the $x_i$ direction in $\R^d$.
Let $x_0\in \R^d$ and let $B=B(x_0,1)$.
For simplicity, we write $\tau$ for $\tau_B$.
We will use $A(x)^{-1}$ to denote the inverse matrix of $A(x)$.

\begin{proposition}\label{P1}
There exist
 positive constants $c_1, c_2$ that depend only on the upper bound of $A(x)$
 and $A (x)^{-1}$ on $B$
such that

\noindent (a) $\E^x [\tau ] \leq c_1$ for all $x\in B$;

\noindent (b) $\E^x [\tau] \geq c_2$ for all $x\in B(x_0, \frac12)$.
\end{proposition}

\proof (a) Let $A_0=\inf\{|A(x)(e_1)|: x\in
 \overline B\}$.
 We know $A_0>0$ because
 $A (x)$ is continuous in $x$ and nondegenerate for
each $x$.
Since the $Z^i$'s are independent one-dimensional symmetric
$\alpha$-stable process, no two of them make a jump at the same time.
So   there exists a positive constant $c_3$ such that
$$
 \P\left(\exists s\leq 1: \Delta Z_s^1>3/A_0 \hbox{ but }
\Delta Z_s^k=0 \hbox{ for } k=2, \cdots, d \right)\geq c_3.
$$
 Suppose there exists  $s\in [0,1]$ such that $\Delta Z_s^1> 3/A_0$,
   $\Delta Z_s^k=0$ for $k=2, \cdots , d$,
and $X_{s-}\in B$. Then by \eqref{RE1}
$$
 |\Delta X_s^1|=
  |\Delta Z_s^1| \, |A(X_{s-})e_1| >3
  $$
 if $X_{s-}\in \overline B$.
 So with probability at least $c_3$,
$X$ will have left $B$ by time
1. Hence if $x\in B$,
$$\P^x (\tau>1)\leq 1-c_3.$$
Let  $\{\theta_t, t>0\}$ denotes the usual shift operators for $X$.
By the Markov property,
\begin{align*}
\P^x(\tau>m+1)&\leq \P^x(\tau>m, \tau\circ\theta_m>1)\\
&=\E^x[\P^{X_m}(\tau>1); \tau>m]\\
&\leq (1-c_3) \P^x(\tau>m).
\end{align*}
By induction,
$$\P^x(\tau>m)\leq (1-c_3)^m,$$
and (a) follows.

(b) Let
$$
\wt Z^i_t:=\sum_{s\leq t} \Delta Z_s^i 1_{(|\Delta Z_s^i|>1)}
\qquad \hbox{and} \qquad \ol Z_t^i :=Z_t^i-\wt Z_t^i.
$$
Note
$$\E[\ol Z^i, \ol Z^i]_t=t\int_{-\beta}^\beta x^2 \frac{c_4}{|x|^{1+\al}} dx
=c_5t\beta^{2-\al}.$$
Let
$\ol X$ solve
$$d\ol X_t=A(\ol X_t) \, d\ol Z_t.$$
 Note that for each $i=1, \cdots, d$, $\ol X^i$ is a purely discontinuous
 square integrable martingale with
 $|\Delta \ol X_t^i|\leq c_6\sum_{j=1}^d |\Delta \ol Z_t^j|$. Hence
$$[\ol X^i, \ol X^i]_t\leq c_7 \sum_{j=1}^d [\ol Z^j, \ol Z^j]_t.$$
First by Chebyshev's inequality and then by Doob's inequality,
\begin{align*}
\P^x\Big(\sup_{s\leq t} |\ol X^i_s-\ol X^i_0|> \frac{1}{4d}\Big)
 &\leq 16 d^2 \, \E \left[ \sup_{s\leq t}|\ol X_s^i-\ol X_0^i|^2 \right]\\
 &\leq 64 d^2 \, \E \left[ (\ol X_t^i-\ol X_0^i)^2 \right]\\
 &=64
d^2\, \E[\ol X^i, \ol X^i]_t\\
&\leq c_8 \sum_{j=1}^d \E[\ol Z^j, \ol Z^j]_t\\
&\leq c_9t.
\end{align*}
Choose $t$ small so that $c_9t\leq 1/4$.

We can choose $t$ smaller if necessary so that
$$\P(\wt Z_s^j\ne 0 \mbox{ for some }s\in [0,t])\leq 1/(4d).$$
So there exists $t$ such that
$\P(\ol Z_s\ne Z_s \mbox{ for some }s\in [0,t])\leq 1/4,$
and it follows that
$$\P(\ol X_s\ne X_s \mbox{ for some }s\in [0,t])\leq 1/4.$$
Therefore with probability at least $1/2$ we have
$\sup_{s\leq t} |X_s-X_0|\leq 1/4$
 and so in particular
$$\P^x(\tau>t)\geq 1/2 \qq \hbox{for } x\in B(x_0, \tfrac12).$$
Consequently,  we have $\E^x \tau \geq t\P^x(\tau\geq t)\geq t/2$
for $x\in B(x, \tfrac12)$.
\qed

\begin{proposition}\label{P2}
 There exist constants $\eta_0>0, p_0\geq 2, $ and $c_1$
 that depend only on the upper bound of $A(x)$
 and $A (x)^{-1}$ on $B$
 such that if
the oscillation of $A$ on $B(x_0,1)$ is less than
$\eta_0$, then
$$\E^x \left[\int_0^\tau 1_C(X_s)\, ds \right]
\leq c_1 |C|^{1/p_0}, \qq x\in B.$$
\end{proposition}

\proof Note that  the process $\{X_t, t\leq \tau\}$ is determined by
the matrix $A$ on $B$ only. Without loss of generality, for this
proof we redefine $A$ for $x\notin B$ so that $A$ is continuous on
$\R^d$ and
$$ \eta:=\sup_{x\in \R^d}\|A(x)-A(x_0)\|= \sup_{x\in B}\|A(x)-A(x_0)\|.
$$
Let $R_\lam$ and $\sL_0$  be the resolvent and infinitesimal
generator of the Levy process $Y_t=Y_0+A(x_0) Z_t$, $\sL$ the
infinitesimal generator of $X$, $S_\lam$ the resolvent of $X$,
and $\sB := \sL -\sL_0$. There exist
$\eta_0>0$ and $p_0\geq 2$ so that the conclusion of Proposition 5.2
of \cite{system} holds, namely, $\norm{\sB R_\lam f}_{p_0}\leq \tfrac14
\norm{f}_{p_0}$.
For $f\in L^{p_0}(\R^d)$, set $h=f-\lam
R_\lam f$. Note that $R_\lam f =R_0h$ and  $\|h\|_{p_0}\leq 2 \|
f\|_{p_0}$. Hence for $\eta<\eta_0$, by \cite[Proposition
5.2]{system}
$$ \norm{\sB R_\lam f}_{p_0} =\norm{\sB R_0 h}_{p_0} \leq
\frac14\|h\|_{p_0}\leq \frac12 \|f\|_{p_0}.
$$
Moreover  by \cite[Proposition 2.2]{system},
$$\norm{R_\lam f}_\infty\leq c_2\norm{f}_{p_0}.$$
It follows from \cite[Proposition 6.1]{system} that
$$S_\lam f=R_\lam \Big(\sum_{i=0}^\infty (\sB R_\lam) ^i\Big)f$$
for $f\in L^{p_0}$ and therefore
\begin{align*}
\norm{S_\lam f}_\infty  =\Bigl\| R_\lam \Big( \sum_{i=0}^\infty
(\sB R_\lam)^i\Big)f\Bigr\|_\infty
 \leq c_2 \Bigl\| \Big(\sum_{i=0}^\infty (\sB R_\lam)^i\Big)f\Bigr\|_{p_0}
 \leq 2c_2 \norm{f}_{p_0}.
\end{align*}
If we apply this to $f=1_C$, where $C\subset B$, then
\begin{equation}\label{E2.1}
\E^x \left[\int_0^\infty e^{-\lam t} 1_C(X_t)\, dt \right]\leq
2c_2|C|^{1/p_0}.
\end{equation}

Let $M=\sup_{x\in B} \E^x \left[ \int_0^\tau 1_C(X_s) \, ds\right]$.
 Clearly $M\leq \sup_{x\in B} \E^x \left[ \tau\right]$,
 which
is finite by Proposition \ref{P1}. By taking $t_1$ sufficiently large,
$$\P^x(\tau\geq t_1)\leq \frac{\sup_{x\in B}\E^x [\tau ]}{t_1}  \leq \tfrac12.$$
We then have
\begin{align*}
\E^x \left[\int_0^\tau 1_C(X_s)\, ds \right]
&\leq \E^x\Big[\int_0^{t_1} 1_C(X_s)\, ds\Big]
+\E^x\Big[ \int_{t_1}^\tau 1_C(X_s)\, ds; \tau\geq t_1\Big]\\
&\leq e^{\lam t_1} S_\lam 1_C(x) +\E^x\left[\E^{X_{t_1}}
\Big[ \int_0^\tau 1_C(X_s)\, ds\Big] ; \tau\geq t_1 \right]\\
&\leq c_3 |C|^{1/p_0} + M\P^x(\tau\geq t_1).
\end{align*}
Taking the supremum over $x$, we have
$$M\leq c_3|C|^{1/p_0} +\tfrac12 M,$$
and our result follows.
\qed

We now prove a support theorem for $X$.
First we prove some lemmas.

\begin{lemma}\label{SPL1}
Let $x_0\in \R^d$, $1\leq k\leq d$,
$v_k=  A(x_0)e_k$,
 $\gamma\in (0,1)$, $t_0>0$, and
$r\in[-1,1]$. There exists $c_1$ depending only on $\gamma$, $t_0$,
$r$,  and the upper bounds and modulus of continuity of $A(\cdot)$
in $B(x_0,2)$ such that
\begin{align}
\P^{x_0}\big(&\mbox{there exists a stopping time }  T\leq t_0 \mbox{ such that } \label{WE1}\\
& \sup_{s<T} |X_s-x_0|<\gamma \hbox{ and } \sup_{T\leq s\leq t_0}
|X_s-(x_0+rv_k)|<\gamma \big)\geq c_1. \nn
\end{align}
\end{lemma}

\proof
Let $\| A\|_\infty:=1 \vee \left( \sum_{i,j=1}^d \sup_{x\in B(x_0, 2)} |A_{ij}(x)|\right)$.
We do the case where $r\geq 0$, the other case being similar. We first suppose
$r\geq \gamma/3$.
Let  $\beta\in (0,r)$ be chosen later,  let
$$\wt Z^i_t=\sum_{s\leq t} \Delta Z_s^i 1_{(|\Delta Z_s^i|>\beta)},
\qq \ol Z_t^i=Z_t^i-\wt Z_t^i,$$
and let $\ol X$ be the solution to
$$d\ol X_s=A(\ol X_{s-})\, d\ol Z_s, \qq \ol X_0=x_0.$$
Choose $\delta<\gamma/
 (6\|A\|_\infty)$
 such that
\begin{equation}\label{e:2.3}
\sup_{i,j} \sup_{|x-x_0|<\delta} |A_{ij}(x)-A_{ij}(x_0)|<\gamma/(12d).
\end{equation}
Let
\begin{align*}
C&= \left\{ \sup_{s\leq t_0} |\ol X_s-\ol X_0|\leq \delta \right\},\\
D&=\{\wt Z_s^i=0 \mbox{ for all }s\leq t_0 \hbox{ and } i\ne k, \
 \wt Z_k \mbox{ has a single jump before time } t_0  \\
&~~~~~~~  \mbox{ and its size is in } [r,r+\delta]\},\\
E&=\{\wt Z_s^i=0 \mbox{ for all }s\leq t_0 \hbox{ and }  i=1, \ldots, d\}.
\end{align*}
As in the proof of Proposition \ref{P1},
$$\E[\ol X^i, \ol X^i]_t\leq c_2\sum_{j=1}^d \E[\ol Z^j, \ol Z^j]_t\leq c_3t\beta^{2-\al},$$
and
 by Chebyshev's inequality and Doob's inequality,
\begin{align*}
\P \left(\sup_{s\leq t_0}|\ol X^i_s-\ol X^i_0|>\delta/\sqrt d \right)
 &\leq \frac{ \E \Big[ \sup_{s\leq t_0}\left( \ol X^i_{s}-\ol X^i_0 \right)^2 \Big]}
 {\delta^2/d}\\
 &\leq  \frac{4  \E \Big[ \left( \ol X^i_{t_0}-\ol X^i_0 \right)^2 \Big]}
 {\delta^2/d}
 \leq \frac{c_4t_0\beta^{2-\al}}{\delta^2}.
\end{align*}
We choose $\beta<r$ so that
\bee\label{Q1}
c_4t_0\beta^{2-\al}\leq
 \delta^2 /(2d),
\eee
 and then
$\P^{x_0}(C)\geq 1/2$.

In order for $\wt Z^k$ to have a single jump before
time $t_0$, and for that jump's size to be in the interval $[r,r+\delta]$,
then by time $t_0$, (a)  $\wt Z^k$ must have no negative jumps; (b) $\wt Z^k$ must have no jumps
whose size lies in $[\beta, r)$; (c) $\wt Z^k$ must have no jumps whose size lies in $(r+\delta,\infty)$;
and (d) $\wt Z^k$ must have precisely one jump whose size lies in the interval $[r,r+\delta]$.
The events described in (a)--(d) are independent and are the probabilities
that Poisson random variables of parameters $c_{5}t_0\beta^{-\al}$,
$c_{5}t_0(\beta^{-\al}-r^{-\al})$, $c_{5}t_0(r+\delta)^{-\al}$, and
$c_{5}t_0(r^{-\al}-(r+\delta)^{-\al})$, respectively, take the  values $0$, $0$, $0$, and $1$,
respectively.
For $j\ne k$, the probability that $\wt Z^j$ does not have a jump
before time $t_0$ is the probability that a Poisson random variable
with parameter $2c_{5}t_0\beta^{-\al}$ is equal to 0.
Since the $\wt Z^j$, $j=1, \cdots, d$, are independent, we thus see that the probability of
$D$ is bounded below by a constant depending on $r,\delta, t_0$  and
$\beta$.
Because the $\ol Z_t$'s are independent of the $\wt Z^j$'s, then $C$
and $D$ are independent. Therefore

\bee\label{WE2}
\P^{x_0}(C\cap D)\geq c_6/2.
\eee
A similar (but slightly easier)  argument shows that
\bee\label{WE3}
\P^{x_0}(C\cap E)\geq c_7.
\eee

If $T$ is the time when $\wt Z^k$ jumps, then $Z_{s-}=\ol Z_{s-}$
for $s\leq T$, and hence $X_{s-}=\ol X_{s-}$ for $s\leq T$. So up to time $T$, $X_s$ does not move more than $\delta$ away from its starting point.
We have
$$\Delta X_T=A(X_{T-}) \Delta Z_T,$$
 so using \eqref{e:2.3} we have that on $C\cap D$,
\begin{align*}
|X_T-(x_0&+rv_k)| \\
&\leq |X_{T-}-x_0|+|\Delta X_T-rv_k|\\
&= |X_{T-}-x_0|+|A(X_{T-})\Delta Z_T-rv_k|\\
&\leq |X_{T-}-x_0|+r|(A(X_{T-})-A(x_0))e_k|
  + |A(X_{T-}) (\Delta Z_T-re_k )|\\
&\leq \delta + rd \gamma/(12d) +\delta
 \|A\|_\infty
<\gamma/2.
\end{align*}

We now apply the strong Markov property at time $T$. By \eqref{WE3},
 $\P^{X_T}
(C\cap E)\geq c_7$ and so
$$
 \P \left(\sup_{T\leq s\leq T+t_0} |X_s-X_T|<\delta \right)\geq c_{8}.
$$
Using the strong Markov property, we have our result with
$c_1=c_7c_{8}/2$.

If $r<\gamma/3$, the argument is easier. In this case we can take
$T$ identically 0, and use \eqref{WE3}. The details are left to the reader.
\qed

\begin{lemma}\label{WL1} Suppose
 $u, v$
 are two vectors in $\R^d$, $\eta\in (0,1)$,  and $p$
is the projection of $v$ onto
 $u$.
 If $|p|\geq \eta|v|$, then
$$|v-p|\leq \sqrt{1-\eta^2} \, |v|.$$
\end{lemma}

\proof
 Note that the vector $v-p$ is orthogonal to the vector $p$.
So by the Pythagorean theorem,
$| v-p|^2= |v|^2 -|p|^2 \leq (1-\eta^2)|v|^2$. \qed

\begin{lemma}\label{WL2} Suppose the entries of $A$ and $A^{-1}$ are bounded by $\Lambda$.
Let $v$ be a vector in $\R^d$,
 $u_k=Ae_k$,
and $p_k$ the projection of $v$ onto $u_k$ for $k=1, \ldots, d$. Then there exists
$\rho\in (0, 1)$ depending only on $\Lambda$ such that for some $k$,
$$|v-p_k|\leq \rho|v|.$$
\end{lemma}

\proof Since the entries of $A^{-1}$ are bounded, then
$|(A^T)^{-1}w|\leq c_1|w|$. Setting $x=(A^T)^{-1}w$, we see
$|A^Tx|\geq c_2|x|$ for any vector $x$.

Let $b_k$ be the projection of $A^Tv$ onto $e_k$. If $|b_k|<(1/d)|A^Tv|$ for all $k$,
then
$$|A^Tv|=\Big|\sum_{k=1}^d b_k\Big|\leq \sum_{k=1}^d |b_k|< |A^Tv|,$$
a contradiction. So for some $k$, $|b_k|\geq (1/d)|A^Tv|\geq c_3|v|$, where
$c_3=c_2/d$.
We then write
$$c_3|v|\leq |b_k|=|v^TAe_k|\leq \frac{c_4}{|Ae_k|} \, |v^TAe_k|
=c_4\frac{|v^Tu_k|}{|u_k|}
 = c_4 |p_k|.
$$
 We now apply Lemma \ref{WL1} with $\eta=c_3/c_4$ and
set $\rho=\sqrt{1-(c_3/c_4)^2}$.
\qed

\begin{lemma}\label{SPL2}
Suppose the entries of $A(x)$ and $A(x)^{-1}$ on $B(x_0, 3)$ are
bounded by $\Lambda$.
Let $t_1>0$,
 $\eps \in (0, 1)$,
$r\in (0,\eps/4)$  and $\gamma>0$. Let
$\psi:[0,t_1]\to \R^d$ be a line segment of length $r$ starting at
$x_0$. Then there exists $c_1>0 $
 that depends only on $\Lambda$, the modulus of continuity of
 $A(x)$ on $B(x_0, 3)$, $t_1$, $\eps$ and $\gamma$
such that
$$
  \P^{x_0} \left( \sup_{s\leq t_1} |X_s-\psi(s)|<
   \eps
  \, \hbox{ and }
  \, |X_{t_1}-\psi(t_1)|<\gamma \right)  \geq c_1.
 $$
\end{lemma}

\proof  Use the bounds on $A$ in $B(x_0,2)$ and Lemma \ref{WL2} to
define $\rho \in (0, 1)$ so that the conclusion of Lemma \ref{WL2}
holds for all matrices $A=A(x)$ with $x\in B(x_0, 2)$. Take
$\gamma\in(0,r\land \rho)$ smaller if necessary so that $\wt
\rho:=\gamma+\rho<1$. Choose $n\geq 2$ large so that $(\wt
\rho)^n<\gamma$.

Let
 $v_0= \psi(t_1)-\psi(t_0)=\psi (t_1)-x_0$, which has length $r$.
 By Lemma \ref{WL2}, there exists $k_0\in \{1, \cdots, d\}$
 such that if $p_0$ is the projection
of $v_0$ onto $A(x_0)e_{k_0}$, then $|v_0-p_0|\leq \rho|v_0|$.
Note $|p_0|\leq |v_0|
 = r$.

Let $D_1$ denote the event that there is a stopping time $T_0\leq t_1/n$
such that
  $|X_s-x_0|<
\gamma^{n+1}$ for $s<T_0$ and $|X_s-(x_0+p_0)|<\gamma^{n+1}$ for $s\in [T_0, \, t_1/n ]$.
By Lemma \ref{SPL1} there exists $c_2>0$ such that $\P^{x_0}(D_1)\geq c_2$.
Note that
on $D_1$,
if $T_0\leq s\leq t_1/n$,
\bee\label{r2E1}
  |\psi(t_1)-X_s|\leq |\psi (t_1)-(x_0+p_0)|+|(x_0+p_0)-X_{t_1/n}|
    \leq \rho r+ \gamma^{n+1} \leq \wt \rho \, r.
\eee
 Taking $s=t_1/n$, we have
$$|\psi(t_1)-X_{t_1}|\leq \wt \rho r.$$
Since $\wt \rho <1$ and
 $|\psi (t_1)-x_0|=r$,
then \eqref{r2E1} shows that
 on $D_1$,
 $$ X_{s}\in B(x_0,2r)\subset B(x_0,\eps/2) \qquad \hbox{if }
  T_0\leq s\leq t_1/n.
  $$
   If $0\leq s<T_0$, then
$|X_s-x_0|<\gamma^{n+1}<r$, and so
   $\{X_s, s\in [0, t_1/n]\}\subset B(x_0, 2r)\subset B(x_0,
  \eps /2)$.

Now let $v_1=\psi(t_1)-X_{t_1/n}$.
When $X_{t_1/n}\in B(x_0, \eps/2)$,  by Lemma \ref{WL2}, there
 exists $k_1$ such that if $p_1$ is the projection
of $v_1$ onto $A(X_{t_1/n})e_{k_1}$, then
$|v_1-p_1|\leq \rho |v_1|$.
 Let $D_2$ be the event that there exists
a stopping time $T_1\in [t_1/n,2t_1/n]$ such that
$|X_s-X_{t_1/n}|<  \gamma^{n+1}$
for $t_1/n\leq s<T_1$ and
$|X_s- (X_{t_1/n}+p_1)| < \gamma^{n+1}$ for $T_1\leq s\leq 2t_1/n$.
Using the Markov property at time $t_1/n$  and applying Lemma \ref{SPL1} again, there exists
(the same) $c_2>0$ such that
$$
 \P^{x_0}(D_2\mid \sF_{t_1/n})\geq c_2
 $$
on the event $\{X_{t_1/n}\in B(x_0,\eps/2)\}$, where $\sF_t$ is the minimal
augmented filtration for $X$. So
$$\P^{x_0}(D_1\cap D_2)\geq c_2\P^{x_0}(D_1)\geq c_2^2.$$
 On the event $D_1\cap D_2$, if $T_1\leq s\leq 2t_1/n$,
\begin{eqnarray*}
 |\psi (t_1)-X_{s}|&\leq& |\psi (t_1)-(X_{t_1/n}+p_1)|+
|(X_{t_1/n}+p_1)-X_{s}| \\
&\leq & \rho|v_1|+\gamma^{n+1} \leq \rho \wt \rho \, r +
\gamma^{n+1} \leq \wt \rho^2 \, r.
\end{eqnarray*}
In particular
$$
   |\psi(t_1)-X_{2t_1/n}|\leq \wt \rho^2 r
   \qquad \hbox{on } D_1\cap D_2.
$$
If $T_1\leq s\leq 2t_1/n$, then $|\psi(t_1)-X_s|<r$ and
$ |\psi(t_1)-
 x_0
 |=r$, and so  $X_s\in B(x_0,2r)\subset B(x_0, \eps/2)$.
 In particular,
 $$ |X_{2t_1/n}-x_0|<\eps /2 \qquad \hbox{on } D_1\cap D_2.
 $$
If $t_1/n\leq s<T_1$, then $|X_s-X_{t_1/n}|<r$ and
$|X_{t_1/n}-x_0|<2r$.
 So on $D_1\cap D_2$,  $X_s\in B(x_0,3r)\subset B(x_0,
  3\eps /4)$.

Let $v_2=\psi(t_1)-X_{2t_1/n}$, and proceed as above
to get events $D_3, \cdots, D_k$. At the $k^{th}$ stage,
 we have
 $$\P^{x_0}(D_k \mid \sF_{(k-1)t_1/n})\geq c_2
 $$
 and so $\P^{x_0}(\cap_{j=1}^k D_j ) \geq c_2^k$.
 On
the event $\cap_{j=1}^k D_j$, if $kt_1/n\leq T_k\leq s\leq
(k+1)t_1/n$, then
$$|\psi(t_1)-X_s|\leq \wt \rho^{k+1} r<r;$$
since $|\psi(t_1)-x_0|=r$, then $X_s\in B(x_0,2r)\subset B(x_0,
\eps/2)$. At the $k^{th}$ stage, on the event $\cap_{j=1}^k D_j$,
$$ | X_{kt_1/n}-x_0|<\eps /2
$$
and
  if $kt_1/n\leq s< T_k$, then
$$|X_s-x_0|\leq |X_s-X_{kt_1/n}|+|X_{kt_1/n}-\psi(t_1)|+|\psi(t_1)-x_0|
< \gamma^{n+1}+2r+r<3r,$$
and so $X_s\in B(x_0,3r)\subset B(x_0,
 3\eps /4)$.

We continue this procedure $n$ times to get events $D_1, \cdots, D_n$
so that on $\cap_{k=1}^n D_k$, we have
 $X_s\in B(x_0, 3\eps /4)$ for $s\leq t_1$,
  $|X_{t_1}-\psi(t_1)|<\gamma$, and $\P^{x_0}(\cap_{k=1}^nD_k)\geq  c_2^n$.
Consequently, on $\cap_{k=1}^n D_k$,
$$ |X_s-\psi (s)|\leq |X_s-x_0|+|x_0-\psi (s) |<3\eps/4 +r < \eps
\qquad \hbox{for } s\in [0, t_1].
$$
This completes the proof of the lemma.
 \qed

\begin{theorem}\label{T3}
Suppose the entries of $A(x)$ and $A (x)^{-1}$ on $B(x_0, 3)$ are
bounded by $\Lambda$.
Let $\vp:[0,t_0]\to \R^d$ be continuous with $\vp(0)=x_0$ and the
image of $\vp$ contained in $B(0,1)$. Let $\eps>0$. There exists
 $c_1>0$ depending on $\Lambda$, the modulus of continuity of
 $A(x)$ on $B(x_0, 3)$,
$\vp, \eps$, and $t_0$ such that
$$
 \P^{x_0} \left(\sup_{s\leq t_0} |X_s-\vp(s)|<\eps \right)>c_1.
$$
\end{theorem}

\proof
We may approximate $\vp$ to within $\eps/2$ by a
polygonal path, so by changing $\eps$ to $\eps/2$, we may without
loss of generality assume $\vp$ is polygonal. Let us now choose $n$ large and
subdivide $[0,t_0]$ into $n$ equal subintervals so that over
each subinterval $[kt_0/n, (k+1)t_0/n]$ the image of
$\vp$ is a line segment of length less than $\eps/4$. We then use
Lemma \ref{SPL2} and the strong Markov property $n$ times to show
 that, with probability at least $c_1>0$,
on each time interval $[kt_0/n, (k+1)t_0/n]$, $X_t$ follows within
$\eps/2$ the line segment from $X_{kt_0/n}$ to $\vp((k+1)t_0/n)$ and
is at most $\eps/(4\sqrt d)$ away from $\vp((k+1)t_0/n)$.
\qed

\begin{corollary}\label{CSP}
Let $\eps, \delta \in (0, 1/4)$. Suppose $Q$ represents either the
unit ball or the unit cube, centered at $x_0\in \R^d$.
Suppose the entries of $A$ and $A^{-1}$ on $Q$ are bounded by
$\Lambda$.
Let $Q'$ be the ball (resp., cube) with radius (resp., side length)
$1-\eps$ with the same center. Let $R$ be a ball (resp., cube) of
radius (resp., side length) $\delta$ contained in $Q'$. Then there
exists $c_1>0$
depending on $\Lambda$, the modulus of continuity of
 $A(x)$ on $Q$, $\eps$ and $\delta$
 such that
$$\P^x(T_R<\tau_Q)\geq c_1, \qq x\in Q'.$$
\end{corollary}

\proof
 Note that the above probability is determined by the
values of the matrix $A(x)$ only on $Q$
 so we can redefine $A(x)$ outside of $Q$ if necessary to make the entries
 of $A$ and $A^{-1}$ on $\R^d$ bounded by $\Lambda$, and the modulus of
 continuity of $A(x)$ on $\R^d$ be the same as that on $Q$.
To prove the corollary, we need only observe that the estimates in
Theorem \ref{T3} can be made to hold uniformly over every line
segment from $x$ to $y$, with $x\in Q'$ and $y$ being the center of
$R$. \qed

A scaling argument shows that
 for $\lam>0$, $\{\wh X_t:=\lam X_{t/\lam^\alpha}, \, t\geq 0\}$
is a process
of the same type as
 $X$.
More precisely, one can show that there exist $d$
 independent one-dimensional symmetric stable processes $\wh Z^j$  of index
$\al$ such that $\wh X$ satisfies
$$d\wh X^i_t=\sum_{j=1}^d \wh A_{ij}(\wh X_t) \, dZ^j_t, \qq \wh X^i_0=
  \lam x_0^i, $$
where $\wh A_{ij}(x)=  A_{ij}(x/\lam)$.
 Note in particular that
 when $\lam \geq 1$,
the oscillation of $\wh A$ will be no more than
the oscillation of $A$. A consequence is that the analogues of
Propositions \ref{P1} and \ref{P2}
and Theorem \ref{T3} hold in balls $B(x_1, r)$ with the same constants
provided $r<1$ (so that $\lam =1/r>1$).

We now have what is needed to prove our main theorem.

\begin{theorem}\label{T4}
 Let $r \in (0, 1]$ and $\gamma >1$.
Suppose $h$ is harmonic in
 $B(x_0, \, \gamma \, r)$
with respect to $X$ and $h$ is bounded  in $\R^d$. There exists
 positive constants $c_1$ and $\beta$ that depend on $\gamma $,
 the upper bound of $A(x)$ and $A(x)^{-1}$ on $B(x_0, \, \gamma\,  r)$,
 and the modulus of continuity of $A(x)$ on $B(x_0, \, \gamma\, r )$ but
 otherwise is independent of $h$ and $r$
 such that
$$|h(x)-h(y)| \leq c_1 \Big(\frac{|x-y|}{r}\Big)^\beta
 \sup_{\R^d} |h(z)|
$$
\end{theorem}

\proof If one examines the proof of Krylov-Safonov carefully (see, e.g., the
presentation in \cite{Ba97}, Theorem V.7.4), one sees that one needs the support theorem and Corollary \ref{CSP},
a result such as Proposition \ref{P2} and
estimates such as Proposition \ref{P1}
and that with these ingredients, one can conclude that if
$Q$ is a cube of side length $r\leq 1$, $A\subset Q\subset B(x_0,r)$,
and $Q'$ is a cube with the same center as $Q$ but side length half as long,
then
\begin{equation}\label{bbb1}
\P^x(T_A<\tau_Q)\geq \vp(|A|/|Q|)  \qq \hbox{for } x\in Q',
\end{equation}
where $\vp$ is a strictly increasing function with $\vp(0)=0$.

Now let $B=B(y,s)$ be a ball contained in
$B(x_0,r)$ and suppose $A\subset B$ with $|A|/|B|\geq 1/3$. Let
$B'=B(y, (1-\eps)s)$, where $\eps$ is chosen so that $|B\setminus B'|/|B|
=1/6$. Then $|A\cap B'|/|B|\geq 1/6$. Cover $B'$ with $N$
equally sized cubes whose interiors are disjoint and each contained
in $B$. We may choose $N$ independent of $s$. For at least one
cube, say, $Q$, we must have $|A\cap B'\cap Q|/|Q|\geq 1/6$. Let
$Q'$ be the cube with the same center as $Q$ but side length half
as long. By the support theorem, if $x\in B(y,s/2)$, there is
probability at least $c_2$ such that
$\P^x(T_{Q'}<\tau_B)\geq c_2$. Applying
(\ref{bbb1}) and the strong Markov property, we have
\begin{equation}\label{bbb2}
\P^x(T_A<\tau_B)\geq c_3>0  \qq \hbox{for } x\in B(y,s/2).
\end{equation}

Applying (\ref{bbb2}) and Proposition \ref{P1}, the result now follows
exactly as the proof in Theorem 4.1 of \cite{BL1}.
(We remark that line 15 on page 386 of \cite{BL1} should read instead
$$(b_{k-1}-a_{k-1})\P^y(\tau_k<T_A)\leq \frac{1}{\gamma}(b_k-a_k)(1-\P^y(T_A<\tau_k)).$$
With suitable modifications to the definition of $\gamma$ and $\rho$, the proof of Theorem 4.1 in
\cite{BL1} is valid.)
\qed

\section{A counterexample to the Harnack  inequality}\label{S4}

We now show that one cannot expect a Harnack inequality to hold, even
when $A(x)\equiv I$, the identity matrix.
We will describe $\eps$ in a moment. Write points in $\R^3$ as $w=(x,y,z)$
and let $w_0=(0, \frac12,0)$.
 Write
$B$ for $B(0,1)$, $\tau$ for $\tau_B$, and let $F_\eps=(-\eps,\eps)^2\subset \R^2$,
$C_\eps=(\R\times F_\eps)\cap B$,
 and $E_\eps=(2,4)\times F_\eps$.            Let $X_t, Y_t $ and $Z_t$ be
independent
one-dimensional symmetric $\alpha$-stable processes and set $W_t=(X_t,Y_t,Z_t)$.
Define $h_\eps(w)=
\P^w(W_{\tau } \in E_\eps)$.
We will show that   $h_\eps(0)/h_\eps(w_0)\to \infty$
as $\eps\to 0$; this
 implies
that a Harnack inequality is not possible.

The L\'evy measure $n(w, d\wt w)$ of $W$ is
$$ n(w, d\wt w)= c \sum_{k=1}^3 |w_k-\wt w_k|^{-1-\alpha}\, d\wt w_k \,
\prod_{j\not= k} \delta_{w_j} (d\wt w_j)
$$
 where $\delta_{a}$ denotes the Dirac measure at the point $a$.
Since all jumps of $W$ are in directions parallel to the coordinate axes, the
only way $W_{\tau}$ can be in $E_\eps$ is if $W_{\tau-}$ is in $C_\eps$.
This is the key observation.

We first get an upper bound on $h_\eps$. It is well known that if $p_t(u,v)$
is the transition density for a one-dimensional symmetric stable
process of index $\al$, then $p_t$ is everywhere strictly positive,
is jointly continuous, $p_t(u,v)=t^{-1/\al}p_1(u/t^{1/\al},
v/t^{1/\al})$, and $p_1(u,v)\sim c_1 |u-v|^{-\al-1}$ for
$|u-v|$ large. An integration
gives
\begin{align*}
 \E^{(y,z)} &\left[\int_0^\infty 1_{(-1,1)^2}(Y_s,Z_s)\, ds \right]\\
& \leq   1 +\int_{1}^\infty \Big(\int_{-1}^1 p_t(y,u)\, du\Big)
 \Big(\int_{-1}^1 p_t(z,v)\, dv\Big)  \, ds  <\infty.
\end{align*}
By scaling,
$$\E^{(y,z)}\left[\int_0^\infty 1_{F_\eps}(Y_s,Z_s)\, ds \right]<c_2\eps^\al.$$
By the L\'evy system
formula (see \cite{BL1} or \cite{CK03}),
\begin{align}
\E^w \left[ \sum_{s\leq t\land \tau} 1_{(W_{s-}\in C_\eps, W_s\in E_\eps)} \right]
&=\E^w \left[\int_0^{t\land \tau} 1_{C_\eps}(W_s) n(W_s, E_\eps)\, ds \right] \nn\\
&\leq c_3 \E^w \left[ \int_0^\infty 1_{C_\eps}(W_s) \, ds \right] \nn\\
 &\leq
c_3 \E^{(y,z)} \left[\int_0^\infty 1_{F_\eps}(Y_s,Z_s)\, ds \right] \nn\\
&\leq c_2c_3 \eps^\al.   \label{aa1}
\end{align}
Letting $t\to \infty$, we obtain
\bee\label{bb1}
h_\eps(w)=\P^w(W_\tau\in E_\eps)\leq c_4\eps^\al.
\eee

Next we get a lower bound on $h_\eps(0)$. Let  $C'_\eps=C_\eps\cap\{ |x|<1/2\}$.
By  the L\'evy system formula we have
\begin{align*}
h_\eps(0)&\geq \E^0 \left[ \sum_{s\leq t\land \tau} 1_{(W_{s-}\in C'_\eps, W_s\in E_\eps)} \right]\\
&=\E^0 \left[ \int_0^{t\land \tau} 1_{C'_\eps}(W_s)  n(W_s, E_\eps)\, ds \right]\\
&\geq c_5 \E^0 \left[ \int_0^{t\land \tau} 1_{C'_\eps}(W_s) \, ds \right].
\end{align*}
Letting $t\to \infty$,
$$h_\eps(0)\geq c_5 \E^0 \left[ \int_0^\tau 1_{C'_\eps}(W_s)\, ds \right].$$
By the  scaling property of $\al$-stable processes, if $\ol V$ is a
one-dimensional symmetric $\alpha$-stable  process starting from 0
killed on exiting $[-1/4,1/4]$, then $\eps^{-1}V_t$  has the same
distribution as $\ol U_{t/\eps^\alpha}$, where $\ol U$ is a
one-dimensional symmetric $\alpha$-stable  process starting from 0
killed on exiting $[-1/(4\eps),1/(4\eps)]$. Hence there is a
positive constant $c_6>0$ such that for every $\eps \in (0, 1)$ and
$t\in (0, \eps^\alpha)$,
$$ \P^0( \ol V_t\in [-\eps, \eps])=\P^x (\ol U_{t/\eps^\alpha} \in [-1, 1])\geq c_{6}.$$
Consequently,
$$\E^0 \left[ \int_0^\infty 1_{C'_\eps}(W_s)\, ds \right]
\geq \E^0 \left[ \int_0^{\eps^\al} 1_{C'_\eps}(\ol W_s)\, ds \right]
\geq c_{7}\eps^\al,$$
where $\ol W$ is the process $W$ killed when any of $X,Y$, or $Z$
exceeds $1/4$ in absolute value. Therefore
\bee\label{bb2}
h_\eps(0)\geq c_{8} \eps^\al.
\eee

Let $G=(-1,1)^2\subset \R^2$, write $\wh w$ for $(y,z)$, and
$\wh W_t=(Y_t,Z_t)$.
By the estimates on  the transition densities,
 we see that
 $$u
 (\wh w):=\E^{\wh w} \left[ \int_0^\infty 1_{G}(\wh W_s) ds \right]
 $$
  is bounded and
\bee\label{aa3}
u(\wh w)
\leq\int_0^\infty
\P^{y}(|Y_s|<1)\P^{z}(|Z_s|<1)\, ds
\to 0
\eee
as $|\wh w|\to \infty$.
Similarly, for $\wh w\in G$,
$$u(\wh w)\geq \int_1^2 \P^y(|Y_s|<1)\P^z(|Z_s|<1)\, ds\geq c_{9}.$$

Now
$u(\wh W_{t\land T_B})$ is a bounded supermartingale, so by optional stopping
$$u(\wh w)\ge\E^{\wh w}[u(\wh W_{T_G}); T_G<\infty]\geq c_{9} \P^w(T_G<\infty),$$
and (\ref{aa3}) then implies that
$\P^{\wh w}(T_G<\infty)\to 0$ as ${\wh w}\to \infty$.
Scaling then shows that
$$
\P^{(1/2,0)}(T_{F_\eps}<\infty)\to 0 \qquad \hbox{as } \eps\to 0,
$$
and hence
\begin{equation}\label{aa2}
\P^{w_0}(T_{C_\eps}<\infty)\to 0 \qquad \hbox{as } \eps\to 0.
\end{equation}
Therefore by \eqref{aa1}-\eqref{bb1},
\begin{align*}
h_\eps(w_0) &=\E^{w_0}[h_\eps(W_{T_{C_\eps}}); T_{C_\eps}<\tau] \\
& \leq c_{10} \eps^\al \P^{w_0}(T_{C_\eps}<\tau) \\
&\leq c_{11}\, h_\eps(0) \,\P^{w_0} (T_{C_\eps}<\infty).
\end{align*}
This and (\ref{aa2}) shows that $h_\eps(0)/h_\eps(w_0)$ can be made as large as we like
by taking $\eps$ small enough
and so a Harnack inequality for $W$ is not possible.

\medskip

\begin{remark}\label{rem1}
{\rm When $\al<1$, we can construct a two-dimensional example
along the same lines.
}
\end{remark}
\bs

\vskip 0.6truein

%

\end{document}